\documentclass[12pt]{article}
\input{xypic}
\usepackage{amsmath}
\usepackage{amssymb}
\usepackage{latexsym}
\usepackage{enumerate}

\newtheorem{theorem}{Theorem}[section]
\newtheorem{lemma}[theorem]{Lemma}

\newtheorem{definition}[theorem]{Definition}
\newtheorem{example}[theorem]{Example}
\newtheorem{proposition}[theorem]{Proposition}
\newtheorem{corollary}[theorem]{Corollary}

\newenvironment{proof}{\noindent{\sc Proof.}}{$\square$}

\newtheorem{remark}[theorem]{Remark}

\numberwithin{equation}{section}

\newcommand{\mono}{\hookrightarrow}

\newcommand{\ra}{\rightarrow}
\newcommand{\la}{\leftarrow}
\newcommand{\epi}{\mbox{$\to$\hspace{-0.35cm}$\to$}}

\def\rmono{\rto|<\hole|<<\ahook}

\def\umono{\ar@{_{(}->}[u]}
\def\uumono{\ar@{_{(}->}[uu]}

\def\lmono{\ar@{_{(}->}[l]}
\def\llmono{\ar@{_{(}->}[ll]}

\def\depi{\dto|>>\tip}

\newcommand{\M}{{\mathcal M}}
\newcommand{\mO}{{\mathcal O}}

\newcommand{\Z}{{\mathbb Z}}

\begin{document}

\title{Fibrewise nullification and the cube theorem
}

\author{
{\sc David Chataur and J\'er\^ome Scherer}
\thanks{The first author was partially supported by DGESIC grant
PB97-0202.}}

\date{\today}
\maketitle

\begin{abstract}
In this paper we explain when it is possible to construct
fibrewise localizations in model categories. For pointed spaces,
the general idea is to decompose the total space of a fibration as
a diagram over the category of simplices of the base and replace
it by the localized diagram. This of course is not possible in an
arbitrary category. We have thus to adapt another construction
which heavily depends on Mather's cube theorem. Working with model
categories in which the cube theorem holds, we characterize
completely those who admit a fibrewise nullification.
\end{abstract}

\maketitle

\section*{Introduction}
Mather's cube theorem states that the top face of a cube of spaces
whose bottom face is a homotopy push-out and all vertical faces
are homotopy pull-backs is again a homotopy push-out
(\cite[Theorem~25]{MR53:6510}). This theorem is one of the very
few occurences of a situation where homotopy limits and colimits
commute. It is actually related to a theorem of Puppe about
commuting fibers and push-outs (\cite{MR51:1808}), and also to
Quillen's Theorem~B in \cite{MR49:2895}. Doeraene's work on
$J$-categories has incorporated the cube theorem as an axiom in
pointed model categories and allowed him to study the
L.S.-category in an abstract setting (\cite{MR94b:55017}). Roughly
speaking a $J$-category is a model category in which the cube
theorem holds. Such a model category is very suitable for studying
the relationship between a localization functor (constructed by
means of certain homotopy colimits) and fibrations.

Recall that a localization functor in a model category $\M$ is any
coaugmented idempotent functor $L: \M \rightarrow \M$. The
coaugmentation is a natural transformation $\eta: Id \rightarrow
L$. We will only deal with nullification functors $P_A$. In this
context the image of $P_A$ is characterized by the property that
$map(A, P_A X) \simeq *$. We are looking for an existence theorem
of fibrewise nullification, i.e. a construction which associates
to any fibration $F \rightarrow E \rightarrow B$ another fibration
together with a natural transformation
\[\xymatrix{
F \rto\dto^{\eta} & E \rto\dto & B \dto\\
P_A F \rto & \bar E \rto & B }
\]
where $E \rightarrow \bar E$ is a $P_A$-equivalence. This is
achieved by imposing the join axiom for the object $A$: We require
the join $X*A$ to be killed by $P_A$, i.e. $P_A (X*A) \simeq *$,
for any object~$X$.

\medskip

For pointed spaces, the most elegant construction of fibrewise
localization is due to E. Dror Farjoun (in \cite[Theorem
F.3]{dror:book}). His idea is to decompose the total space of a
fibration as a diagram over the category of simplices of the base
and replace it by the corresponding localized diagram. In certain
particular settings, some authors used other constructions (P. May
\cite{MR81f:55005}, W.~Dwyer, H.~Miller, and J.~Neisendorfer in
\cite{MR90i:55034} for completions, C.~Casacuberta and
A.~Descheemaker in \cite{CD} in the category of groups), but none
of these can be adapted in model categories. We prove the
following:

\medskip

\noindent{\bf Theorem \ref{fibrewise}}\,
{\it Let $\M$ be a model category which is pointed, left proper,
cellular and in which the cube and the join axiom hold. Then the
nullification functor $P_A$ admits a fibrewise version.}

\medskip

This condition is actually necessary and we characterize
completely the model categories for which fibrewise nullifications
exist. This is closely related to the property of preserving
products: A nullification functor $P_A$ preserves (finite)
products if $P_A(X \times Y) \simeq P_A X \times P_A Y$.

\medskip

\noindent{\bf Theorem \ref{charac}}\,
{\it Let $\M$ be a model category which is pointed, left proper,
cellular and in which the cube axiom holds. Then the following
conditions are equivalent:
\begin{itemize}
\item[(i)]{} The nullification functor $P_A$ admits a fibrewise
version.
\item[(ii)]{} The nullification functor $P_A$ preserves
finite products.
\item[(iii)]{} The canonical projection $X \times
A \rightarrow X$ is a $P_A$-equivalence for any $X \in \M$.
\item[(iv)]{} The join axiom for $A$ is satisfied.
\end{itemize}
}

\medskip

We show in the last part of the paper that the category of
algebras over an admissible operad satisfies the cube axiom.
Therefore the plus-construction developed in \cite{CRS} has a
fibrewise analogue. Let us only say that the plus-construction
performed on a $\mO$-algebra $B$ kills the maximal $\mO$-perfect
ideal in $\pi_0 B$ and preserves Quillen homology. As a direct
consequence we get the following result which is classical for
spaces.

\medskip

\noindent{\bf Theorem \ref{acyclic}}\,
{\it Let $\mO-alg$ be the category of algebras over an admissible operad $\mO$.
For any $\mO$-algebra $B$, denote by $B \rightarrow B^+$ the plus construction.
The homotopy fiber $AB = Fib(B \rightarrow B^+)$ is then acyclic with respect
to Quillen homology.}

\medskip

\noindent {\bf Acknowledgements}. We would like to thank Gustavo
Granja and Sophie Reinberg for helpful comments.

\section{The cube axiom}
We work in a model category $\M$ which is {\it pointed}, i.e. the
terminal object coincides with the initial one and is denoted by
$*$. In such a category the homotopy fiber $Fib(p)$ of a map $p:E
\rightarrow B$ is defined as the homotopy pull-back of the diagram
$* \rightarrow B \leftarrow E$. We also assume the category is
{\it left proper}, meaning that the push-out of a weak equivalence
along a cofibration is again a weak equivalence. Finally we
require $\M$ to be {\it cellular} as defined in
\cite[Definition~14.1.1]{hirschhorn:unpub}. Basically the small
object argument applies in a cellular model category, as one has
$I$-cells which replace the usual spheres. There exists a cardinal
$\kappa$ such that any morphism from an $I$-cell to a telescope of
length $\lambda \geq \kappa$ factorizes through an object of this
telescope. Moreover every object has a cofibrant replacement by an
$I$-cell complex by \cite[Theorem~13.3.7]{hirschhorn:unpub}.
Localization functors exist in this setting, see
\cite[Theorem~4.1.1]{hirschhorn:unpub}, but in general we do not
know if it is possible to localize fibrewise in any (pointed, left
proper, cellular) model category. We will thus work in model
categories satisfying an extra-condition.

\begin{definition}
\label{cube}
{\rm A model category $\M$ satisfies the {\it cube axiom} if for
every commutative cubical diagram in $\M$ in which the bottom face
is a homotopy push-out square and all vertical faces are homotopy
pull-back squares, then the top face is a homotopy push-out square
as well.}
\end{definition}

M.~Mather proved the cube Theorem for spaces in
\cite[Theorem~25]{MR53:6510} and J.-P.~Doeraene introduced it as
an axiom for model categories. His paper \cite{MR94b:55017}
contains a very useful appendix with several examples of model
categories satisfying this rather strong axiom.

\begin{example}
\label{stablecube}
{\rm Any stable model category satisfies the cube axiom. Indeed
homotopy push-outs coincide with homotopy pull-backs, so that this
axiom is a tautology. On the other hand the category of groups
does not satisfy the cube axiom. Let us give an easy
counter-example by considering the push-out of $(\Z \leftarrow *
\rightarrow \Z)$, which is a free group on two generators $a$ and
$b$. The pull-back along the inclusion $\Z <ab> \mono \Z<a>*
\Z<b>$ is obviously not a push-out diagram. However fibrewise
localizations exist in the category of groups as shown by the
recent work of Casacuberta and Descheemaker \cite{CD}.}
\end{example}

The following proposition claims that under very special
circumstances the push-out of the fibers coincides with the fiber
of the push-outs. In the category of spaces this is originally due
to V.~Puppe, see~\cite{MR51:1808}. The close link between the cube
Theorem and Puppe's theorem was already well-known to M.~Mather
and M.~Walker, as can be seen in \cite{MR82a:55008}.

\begin{proposition}
\label{Puppe}
Let $\M$ be a pointed model category in which the cube axiom
holds. Consider natural transformations between push-out diagrams:
\[\xymatrix{
F\dto_{j}\ar @{}[r]|(0.46)= &
*{hocolim\hspace{2mm}\big(\hspace{-20pt}} & F_{1}\dto_{j_{1}}
&F_{0}\lto\rto\dto_{j_{0}} & F_{2}\dto^{j_{2}} & *{\hspace{-20pt}\big)}\\
E\dto_{p}\ar @{}[r]|(0.46)= &
*{hocolim\hspace{2mm}\big(\hspace{-20pt}} & E_{1}\dto_{p_{1}}
&E_{0}\lto\rto\dto_{p_{0}} & E_{2}\dto^{p_{2}} & *{\hspace{-20pt}\big)}\\
B \ar @{}[r]|(0.46)=& *{hocolim\hspace{2mm}\big(\hspace{-20pt}} &
B & B \ar @{=}[l] \ar @{=}[r] & B & *{\hspace{-20pt}\big)} }\]
Assume that $F_i = Fib(p_i)$ for any $0 \leq i \leq 2$. Then $F =
Fib(p)$.
 \end{proposition}

\begin{proof}
Denote by $k: G \rightarrow E$ the homotopy fiber of $p$. We show
that $G$ and $F$ are weakly equivalent. Let us construct a cube by
pulling-back $E_i \rightarrow E$ along $k$. The bottom face
consists thus in the middle row of the above diagram and the top
face consists in the homotopy pull-backs of $E_i \rightarrow E
\leftarrow G$, which are the same as the homotopy pull-backs of
$E_i \rightarrow B \leftarrow *$, i.e. $F_i$. The cube axiom now
states that the top face is a homotopy push-out and we are done.
\end{proof}

\medskip

This result will be the main tool in constructing fiberwise
localization in $\M$. In his paper \cite{MR94b:55017} on
L.S.-category, J.-P.~Doeraene used the cube axiom in a very
similar fashion to study fiberwise joins. Indeed Ganea's
characterization of the L.S.-category uses iterated fibers of
push-outs over a fixed base space. The same ideas have also been
used in \cite{MR96i:55030}.

\begin{lemma}
\label{product}
Let $\M$ be a model category in which the cube axiom holds. Let
$D$ be the homotopy push-out in $\M$ of the diagram $A \leftarrow
B \rightarrow C$. Then, for any object $X \in \M$, $X \times D$ is
the homotopy push-out of the diagram $X \times A \leftarrow X
\times B \rightarrow X \times C$.
\end{lemma}

\begin{proof}
It suffices to consider the cube obtained by pulling back the
mentionned push-out square along the canonical projection $X
\times D \rightarrow D$.
\end{proof}

\section{The join}
We check here that we can use all the classical facts about the
join in any model category and introduce the join axiom. Most
proofs here are not new, but probably folklore. Recall that the
join $A*B$ of two objects $A, B \in \M$ is the homotopy push-out
of $A \stackrel{p_1}{\longleftarrow} A \times B
\stackrel{p_2}{\longrightarrow} B$. First notice that the induced
maps $A \ra A*B$ and $B \ra A*B$ are trivial. Indeed the map $A
\ra A*B$ can be seen as the composite $A
\stackrel{i_1}{\longrightarrow} A \times B
\stackrel{p_1}{\longrightarrow} A \ra A*B$ which by definition
coincides with the obviously trivial map $A
\stackrel{i_1}{\longrightarrow} A \times B
\stackrel{p_2}{\longrightarrow} B \ra A*B$.

\begin{lemma}
\label{join}
For any objects $A, B \in \M$, we have $A*B \simeq \Sigma(A \wedge B)$.
\end{lemma}

\begin{proof}
We use a ``classical" Fubini argument (homotopy colimit commute
with itself, cf. for example \cite[Theorem~24.9]{ChSc}). Let $P$
be the homotopy push-out of $A \la A \vee B \ra A \times B$ and
consider first the commutative diagram
\[\xymatrix{
A & A \vee B \lto \rto & B\\
A \ar @{=}[u] \ar @{=}[d] & A \vee B \ar @{=}[u] \dto \rto \lto & A \times B \ar @{=}[d] \uto\\
A & A \times B \lto \ar @{=}[r] & A \times B }
\]
Its homotopy colimit can be computed in two different ways. By
taking first vertical homotopy push-outs and next the resulting
horizontal homotopy push-out one gets $A*B$. By taking first
horizontal homotopy push-outs one gets the homotopy cofiber of $P
\ra A$. Consider finally the commutative diagram
\[\xymatrix{
\ast & \ast \ar @{=}[r] \ar @{=}[l] & \ast\\
A \uto \ar @{=}[d] & A \vee B \uto \dto \rto \lto & A \times B \ar @{=}[d] \uto\\
A & A \times B \lto \ar @{=}[r] & A \times B }\]
\noindent
 The same process as above shows that $Cof(P \ra A)$ is homotopy equivalent
to $\Sigma(A \wedge B)$.
\end{proof}

\begin{lemma}
\label{suspension}
For any objects $A, B \in \M$, we have $\Sigma A \wedge B \simeq
\Sigma(A \wedge B)$.
\end{lemma}

\begin{proof}
Apply again the Fubini commutation rule to the following diagram
\[\xymatrix{
\ast & \ast \ar @{=}[r] \ar @{=}[l] & \ast\\
B \uto \ar @{=}[d] & A \vee B \uto \dto \rto \lto & B \ar @{=}[d] \uto\\
B & A \times B \lto \rto & B }\]
\noindent
where one uses
Lemma~\ref{product} to identify the push-out of the bottom line.
\end{proof}

\medskip

For a fibration $F \ra E \epi B$, the {\it holonomy action} is the map
$m: \Omega B \times F \ra F$ induced on the pull-backs by the natural transformation
from $\Omega B \ra * \la F$ to $PB \epi B \la E$.

\begin{corollary}
\label{holonomy}
For any fibration $F \ra E \epi B$, the homotopy push-out of
$\Omega B \la \Omega B \times F \stackrel{m}{\longrightarrow} F$ is
weakly equivalent to $\Omega B * F$.
\end{corollary}

\begin{proof}
Copy the proof above to compare this homotopy push-out to $\Sigma(\Omega B \wedge F)$.
\end{proof}

\medskip

When working with a nullification functor $P_A$ for some object $A
\in \M$, we say that $X$ is $A$-{\it acyclic} or {\it killed} by
$A$ if $P_A X \simeq *$. By universality this is equivalent to
$map(X, Z) \simeq *$ for any $A$-local object $Z$, or even better
to the fact that any morphism $X \ra Z$ to an $A$-local object is
homotopically trivial.

\begin{definition}
\label{axiom}
{\rm A cellular model category $\M$ satisfies the {\it join axiom}
for the nullification functor $P_A$ if the join of $A$ with any
$I$-cell is $A$-acyclic. }
\end{definition}

\begin{example}
\label{stablejoin}
{\rm Any stable model category satisfies trivially the join axiom,
as push-outs coincide with pull-backs. In such a category the join
is always trivial. The category of groups satisfies the join axiom
for a similar reason (but we saw in Example~\ref{stablecube} that
the cube axiom does not hold).}
\end{example}

\begin{proposition}
\label{acyclicjoin}
Let $\M$ be a cellular model category in which the join axiom and
the cube axiom hold. Then $\Sigma^i A * Z$ is $A$-acyclic for any
$i \geq 0$ and any object $Z$.
\end{proposition}

\begin{proof}
The join is a homotopy colimit and thus commutes with other
homotopy colimits. Since any object in $\M$ has a cofibrant
approximation which can be constructed as a telescope by attaching
$I$-cells, the lemma will be proven if we show that $\Sigma^i A
* Z$ is acyclic for any $I$-cell $Z$. By assumption we know that
$A * Z$ is acyclic and we conclude by Lemma~\ref{suspension} since
$\Sigma^i A * Z \simeq \Sigma^i(A*Z)$ is $P_A$-acyclic.
\end{proof}

\begin{remark}
\label{restricted}
{\rm Given a family $S$ of $I$-cells, we say $\M$ satisfies the
{\it restricted join axiom} if the join of $A$ with any $I$-cell
in $S$ is $A$-acyclic. One refines then the above proposition to
cellular model categories in which the restricted join axiom
holds. Here $\Sigma^i A * Z$ is $A$-acyclic for any $i \geq 0$ and
any $S$-cellular object $Z$, i.e. any object weakly equivalent to
one which can be built by attaching only $I$-cells in $S$. }
\end{remark}

\section{Fibrewise nullification}
Let $A$ be any object in $\M$. Recall that it is always possible
to construct mapping spaces up to homotopy in $\M$ eventhough we
do not assume $\M$ is a simplicial model category (see
\cite{ChSc}). Thus we can define an object $Z \in \M$ to be
$A$-local if there is a weak equivalences $map(A, Z) \simeq *$. A
map $g: X \rightarrow Y$ is a $P_A$-equivalence if it induces a
weak equivalences on mapping spaces $g^*: map(Y, Z) \rightarrow
map(X, Z)$ for any $A$-local object $Z$. Hirschhorn shows that
there exists a coaugmented functor $P_A: \M \rightarrow \M$ such
that the coaugmentation $\eta: X \rightarrow P_A X$ is a
$P_A$-equivalence to an $A$-local object. This functor is called
nullification or periodization.

The nullification $X \ra P_A X$ can be constructed up to homotopy
by imitating the topological construction~2.8 in
\cite{bousfield:local}. One must iterate (possibly transfinitely,
for a cardinal given by the smallness of any cofibrant object in
$\M$, see \cite[Theorem~14.4.4]{hirschhorn:unpub}) the process of
gluing $A$-cells, i.e. take the homotopy cofiber of a map
$\Sigma^i A \ra X$. We assume throughout this section that the
model category $\M$ satisfies both the join axiom and the cube
axiom.

\medskip

Let us explain now how to adapt the fibrewise construction
\cite[F.7]{dror:book} in a model category. The following lemma is
the step we will iterate on and on so as to construct the space
$\bar E$ (in Theorem~\ref{fibrewise}).

\begin{proposition}
\label{step}
Consider a commutative diagram
\[\xymatrix{
F \dto_{j} \rmono^\eta & P_A F \rto\dto_{} & F_{1}\dto^{j_{2}}\\
E \depi_{p} \rmono & E'\rto^{\simeq}\dto_{p'} & E_{1} \depi^{p_{1}}\\
B \ar @{=}[r] & B \ar @{=}[r] & B }\] where the left column is a
fibration sequence, the upper left square is a homotopy push-out
square, $p':E' \ra B$ is the unique map extending $p$ such that
the composite $P_A F \ra E' \ra B$ is trivial, $p_1$ is a
fibration, and $F_1$ is the homotopy fiber of $p_1$. Then the
composites $E \ra E' \ra E_1$ and $F \ra P_A F \ra F_1$ are both
$P_A$-equivalences.
\end{proposition}

\begin{proof}
We can assume that the map $\eta: F \mono P_A F$ is a cofibration
as indicated in the diagram, so that $E'$ is obtained as a
push-out, not only a homotopy push-out. Since $\eta$ is a
$P_A$-equivalence, so is its push-out along $j$ by left properness
(see \cite[Proposition~3.5.4]{hirschhorn:unpub}). To prove that $F
\ra F_1$ is a $P_A$-equivalence, it suffices to analyze the map
$P_A F \ra F_1$. We use Puppe's Proposition~\ref{Puppe} to compute
$F_1$ as homotopy push-out of the homotopy fibers of $P_A F \la F
\ra E$ over the fixed base~$B$. This yields the diagram $P_A F
\times \Omega B \la F \times \Omega B \ra F$ whose homotopy
push-out is $F_1$. We investigate more closely the map $F \ra F_1$
by decomposing the map $F \ra P_A F$ into several steps obtained
by gluing $A$-cells.

Consider a cofibration of the form $\Sigma^i A
\stackrel{f}{\longrightarrow} F \ra C_f$. Let $E_f$ be the
homotopy push-out of $C_f \la F \ra E$ and compute as above the
homotopy fiber $F_f$ of $E_f \ra B$. It is weakly equivalent to
the homotopy push-out of $C_f \times \Omega B \la F \times \Omega
B \ra F$. Hence $F_f$ is also weakly equivalent to the homotopy
push-out of $\Omega B \la \Sigma^i A \times \Omega B \ra F$, using
the definition of $C_f$. Decompose this push-out as follows
\[\xymatrix{
\Sigma^i A \times \Omega B \rto \dto & \Sigma^i A \rto \dto & F \dto \\
\Omega B \rto & \Sigma^i A \ast \Omega B \rto & F_f
}\]
The right-hand square must be a homotopy push-out square as well. But both
$\Sigma^i A$ and $\Sigma^i A \ast \Omega B$ are $A$-acyclic (by Proposition~\ref{acyclicjoin}),
so that the map $\Sigma^i A \ra \Sigma^i A \ast \Omega B$ is a $P_A$-equivalence.
Thus so is $F \ra F_f$ by left properness. Iterating this process of gluing $A$-cells
shows that $F \ra F_1$ is a telescope of $P_A$-equivalences, hence a $P_A$-equivalence.
\end{proof}

\begin{remark}
\label{commutation}
{\rm In the category of spaces it is of course true that $\Omega B
\times F \ra \Omega B \times P_A F$ is a $P_A$-equivalence,
because localization commutes with finite products. In general we
will see in Theorem~\ref{charac} that the join axiom is actually
equivalent to the commutation of $P_A$ with products. With the
restricted join axiom we would have to impose the additional
restriction on $B$ that $\Omega B$ be $S$-cellular.}
\end{remark}

\begin{theorem}
\label{fibrewise}
Let $\M$ be a model category which is pointed, left proper,
cellular and in which the cube axiom and the join axiom hold. Let
$P_A: \M \rightarrow \M$ be a nullification functor. Then there
exists a fibrewise nullification, i.e. a construction which
associates to any fibration $F \rightarrow E \rightarrow B$
another fibration together with a natural transformation
\[\xymatrix{
F \rto\dto^{\eta} & E \rto\dto & B \dto\\
P_A F \rto & \bar E \rto & B
}\]
where $E \rightarrow \bar E$ is a $P_A$-equivalence.
\end{theorem}

\begin{proof}
We construct first by the method provided in Lemma~\ref{step} a
natural transformation to the fibration $F_1 \ra E_1 \ra B$. We
iterate then this construction and get a fibration $\bar F \ra
\bar E \ra B$ where $\bar F = hocolim(F \ra F_1 \ra F_2 \ra
\dots)$ and $\bar E = hocolim(E \ra E_1 \ra E_2 \ra \dots)$. All
maps in these telescopes are $P_A$-equivalences by the lemma,
hence so are $E \ra \bar E$ and $F \ra \bar F$. Moreover any map
$F_n \ra F_{n+1}$ factorizes as $F_n \ra P_A F_n \simeq P_A F \ra
F_{n+1}$ so that $\bar F \simeq P_A F$. We obtain thus the desired
fibration $P_A F \ra \bar E \ra B$.
\end{proof}

\medskip

Define $\bar P_A X = Fib(X \ra P_A X)$, the fiber of the nullification. As in
the case of spaces we get:

\begin{corollary}
\label{A-null}
For any object $X$ in $\M$ we have $P_A \bar P_A X \simeq *$.
\end{corollary}

\begin{proof}
Apply the fiberwise localization to the fibration $\bar P_A X \ra
X \ra P_A X$. This yields a fibration $P_A \bar P_A X \ra \bar X
\ra P_A X$ in which the base and the fiber are $A$-local.
Therefore $\bar X$ is $A$-local as well. But then $\bar X \simeq
P_A X$ and so $P_A \bar P_A X \simeq *$.
\end{proof}

We end this section with a complete characterization of the model
categories which admit fibrewise nullifications.

\begin{theorem}
\label{charac}
Let $\M$ be a model category which is pointed, left proper,
cellular and in which the cube axiom holds. Then the following
conditions are equivalent:
\begin{itemize}
\item[(i)]{} The nullification functor $P_A$ admits a fibrewise
version.
\item[(ii)]{} The nullification functor $P_A$ preserves
finite products.
\item[(iii)]{} The canonical projection $X \times
A \rightarrow X$ is a $P_A$-equivalence for any $X \in
\M$.
\item[(iv)]{} The join axiom for $A$ is satisfied.
\end{itemize}
\end{theorem}

\begin{proof}
We prove first that (i) implies (ii). Consider the trivial
fibration $X \rightarrow X \times Y \rightarrow Y$ and apply the
fibrewise nullification to get a new fibration $P_A X \rightarrow
E \rightarrow Y$. The inclusion of the fiber admits a retraction
$E \rightarrow P_A X$, i.e. $E \simeq P_A X \times Y$. Applying
once again the fibrewise nullification to $Y \rightarrow Y \times
P_A X \rightarrow P_A X$, we see that the map $X \times Y
\rightarrow P_A X \times P_A Y$ is a $P_A$-equivalence. As a
product of local objects is local, this means precisely that
$P_A(X \times Y) \simeq P_A X \times P_A Y$.

Property (iii) is a particular case of (ii). We show now that
(iii) implies (iv). If the canonical projection $X \times A
\rightarrow X$ is a $P_A$-equivalence, the push-out of it along
the other projection yields another $P_A$-equivalence, namely $A
\rightarrow X*A$. Therefore the join $X*A$ is $P_A$-acyclic.
Finally (iv) implies (i) as shown in Theorem~\ref{fibrewise}.
\end{proof}

\medskip

The construction we propose for fibrewise nullification does not
translate to the setting of general localization functors. We do
not know if the cube and join axioms are sufficient conditions for
the existence of fibrewise localizations.

\section{Algebras over an operad}
In this section we provide the motivating example for which this
theory has been developped. For a fixed field $k$, we work with
$\Z$-graded differential $k$-vector spaces ($k$-dgm) and consider
the category of algebras in $k$-dgm over an admissible operad.
This is indeed  a pointed, left proper and cellular category. Weak
equivalences are quasi-isomorphisms and fibrations are
epimorphisms.

We do not know if the join axiom holds in full generality for any
object $A$. It does so however when $A$ is acyclic with respect to
Quillen homology, which is the case we are most interested in, or
when $A$ is a free algebra. We check that the cube axiom always
holds, following the strategy of \cite[Proposition
A.15]{MR94b:55017}, which guarantees the existence of fibrewise
versions of the plus-construction and Postnikov sections. In the
case of $\mathbb N$-graded $\mathcal O$-algebras (the case
$\mathcal O=\mathcal As$ is treated by Doeraene) one has to
restrict to a particular set of fibrations (the so-called
$J$-maps), because they must be surjective in each degree in order
to compute pull-backs. In our context all fibrations are
epimorphisms, so that the cube axiom holds in full generality.

\begin{theorem}
\label{operadiccube}
The cube axiom holds in the category of $\mO$-algebras.
\end{theorem}

\begin{proof}
Let us briefly recall the key steps in Doeraene's strategy. We
consider a push-out square of $\mO$-algebras (along a generic
cofibration $B \mono B \coprod \mO(V)$):
\[\xymatrix{
B \rmono\dto & C=B \coprod \mO(V) \dto \\
A \rmono & D=A \coprod \mO(V)}
\]
We need to compute the pull-back of this square along a fibration
$p: E \epi D$ (which is hence an epimorphism of chain complexes).
We have thus the following isomorphism of chain complexes:
$$E\cong A \coprod \mO(V)\oplus ker(p).$$
This allows to compute the successive pull-backs $A\times_D E$,
$C\times_D E$, and $B\times_D E$. In order to construct the
homotopy push-out $P$ of these pull-backs (which must coincide
with $E$) we factorize the morphism $B\times_{D}E \rightarrow
C\times_{D}E$ as
$$B\times_{D}E\hookrightarrow (B\oplus ker{p})\coprod\mO(V\oplus W)\stackrel{\sim}{\twoheadrightarrow} C\times_{D}E$$
Thus $P$ is identified with $(A\oplus ker(p))\coprod \mO(V\oplus
W)$, which allows us to build finally a quasi-isomorphism to $E$.
\end{proof}

The plus-construction for an $\mO$-algebra is a nullification with
respect to a universal acyclic algebra $\mathcal U$. We refer to
\cite{CRS} for an explicit construction and nice applications.

\begin{proposition}
\label{operadicjoin}
The join axiom holds for any acyclic $\mO$-algebra $A$. It holds
in particular for the universal acyclic algebra $\mathcal U$
constructed in~\cite{CRS}, so that the fibrewise plus-construction
exists.
\end{proposition}

\begin{proof}
The join $A*X$ is weakly equivalent to $\Sigma A \wedge X$ by
Lemmas~\ref{join} and \ref{suspension}. Since $\Sigma A$ is
0-connected and acyclic, it is trivial by the Hurewicz Theorem
\cite[Theorem 1.1]{CRS}. Thus $A*X \simeq *$ is always
$P_A$-acyclic.
\end{proof}

We consider next the case of Postnikov sections $P_{\mO(x)}$,
where $x$ is a generator of arbitrary degree $n \in \Z$. Because
$[\mO(x), X] \cong \pi_n X$ for any $\mO$-algebra $X$, the
nullification functor $P_{\mO(x)}$ is really a Postnikov section,
i.e. $P_{\mO(x)} X \simeq X[n-1]$. Let us also recall that
$\pi_n(X \times Y) \cong \pi_n X \times \pi_n Y$.

\begin{proposition}
\label{spherejoin}
The join axiom holds for any free $\mO$-algebra $\mO(x)$ on one
generator of degree $n \in \Z$. Therefore fibrewise Postnikov
sections exist.
\end{proposition}

\begin{proof}
By Theorem~\ref{charac} we might as well check that the map $X
\times \mO(x) \rightarrow X$ is a $P_{\mO(x)}$-equivalence for any
$\mO$-algebra $X$. Clearly the $n$-th Postnikov section of the
product $X \times \mO(x)$ is equivalent to $X[n-1]$ and we are
done.
\end{proof}

\medskip

Our final result is a particular case of Corollary~\ref{A-null}. A
direct proof (without fibrewise techniques) seems out of reach.

\begin{theorem}
\label{acyclic}
Let $\mO-alg$ be the category of algebras over an admissible
operad $\mO$. For any $\mO$-algebra $B$, denote by $B \rightarrow
B^+$ the plus construction. The homotopy fiber $AB = Fib(B
\rightarrow B^+)$ is then acyclic with respect to Quillen
homology. \hfill{$\square$}
\end{theorem}

\bibliographystyle{alpha}\label{bibliography}
\bibliography{bibho}

\begin{thebibliography}{DMN89}

\bibitem[Bou94]{bousfield:local}
A.~K. Bousfield.
\newblock Localization and periodicity in unstable homotopy theory.
\newblock {\em J. Amer. Math. Soc.}, 7(4):831--873, 1994.

\bibitem[CD02]{CD}
Carles Casacuberta and An~Descheemaker.
\newblock Relative group compeltions.
\newblock {\em preprint}, 2002.

\bibitem[CRS03]{CRS}
David Chataur, Jos\'e Rodr\'{\i}guez, and J\'er\^ome Scherer.
\newblock Plus-construction of algebras over an operad, {H}ochschild and cyclic
  homologies up to homotopy.
\newblock {\em preprint}, 2003.

\bibitem[CS02]{ChSc}
Wojciech Chach\'olski and J\'er\^ome Scherer.
\newblock Homotopy theory of diagrams.
\newblock {\em Mem. Amer. Math. Soc.}, 155(736):ix+90, 2002.

\bibitem[DF96]{dror:book}
E.~Dror~Farjoun.
\newblock {\em Cellular spaces, null spaces and homotopy localization}, volume
  1622 of {\em Lecture Notes in Mathematics}.
\newblock Springer-Verlag, Berlin, 1996.

\bibitem[DMN89]{MR90i:55034}
William Dwyer, Haynes Miller, and Joseph Neisendorfer.
\newblock Fibrewise completion and unstable {A}dams spectral sequences.
\newblock {\em Israel J. Math.}, 66(1-3):160--178, 1989.

\bibitem[Doe93]{MR94b:55017}
Jean-Paul Doeraene.
\newblock L.{S}.-category in a model category.
\newblock {\em J. Pure Appl. Algebra}, 84(3):215--261, 1993.

\bibitem[DT95]{MR96i:55030}
J.-P. Doeraene and D.~Tanr{\'e}.
\newblock Axiome du cube et foncteurs de {Q}uillen.
\newblock {\em Ann. Inst. Fourier (Grenoble)}, 45(4):1061--1077, 1995.

\bibitem[Hir]{hirschhorn:unpub}
P.~S. Hirschhorn.
\newblock Localization of model categories.
\newblock Unpublished, available at P. Hirschhorn's homepage: {\tt
  http://www-math.mit.edu/$\sim$psh/}.

\bibitem[Mat76]{MR53:6510}
Michael Mather.
\newblock Pull-backs in homotopy theory.
\newblock {\em Canad. J. Math.}, 28(2):225--263, 1976.

\bibitem[May80]{MR81f:55005}
J.~P. May.
\newblock Fibrewise localization and completion.
\newblock {\em Trans. Amer. Math. Soc.}, 258(1):127--146, 1980.

\bibitem[MW80]{MR82a:55008}
Michael Mather and Marshall Walker.
\newblock Commuting homotopy limits and colimits.
\newblock {\em Math. Z.}, 175(1):77--80, 1980.

\bibitem[Pup74]{MR51:1808}
Volker Puppe.
\newblock A remark on ``homotopy fibrations''.
\newblock {\em Manuscripta Math.}, 12:113--120, 1974.

\bibitem[Qui73]{MR49:2895}
Daniel Quillen.
\newblock Higher algebraic ${K}$-theory. {I}.
\newblock In {\em Algebraic $K$-theory, I: Higher $K$-theories (Proc. Conf.,
  Battelle Memorial Inst., Seattle, Wash., 1972)}, pages 85--147. Lecture Notes
  in Math., Vol. 341. Springer, Berlin, 1973.

\end{thebibliography}

\bigskip

\bigskip\noindent
David Chataur

\noindent Centre de Recerca Matem\`atica, E--08193 Bellaterra \\
email: {\tt chataur@crm.es}

\bigskip\noindent
J\'er\^ome Scherer

\noindent Departament de Matem\`atiques, Universitat Aut\'onoma de
Barcelona, E--08193 Bellaterra \\
e-mail: {\tt jscherer@mat.uab.es}

\end{document}